\documentclass[a4paper,12pt]{amsart}
\usepackage{amssymb}
\usepackage{amsmath}
\usepackage{ifthen}
\nonstopmode \numberwithin{equation}{section}
\newtheorem{thm}[equation]{Theorem}
\newtheorem{cor}[equation]{Corollary}
\newtheorem{lem}[equation]{Lemma}
\newtheorem{prop}[equation]{Proposition}

\newtheorem{conj}[equation]{Conjecture}
\theoremstyle{definition}
\newtheorem{defn}{Definition}[section]
\newtheorem{example}{Example}[section]
\newtheorem{prob}[equation]{Problem}
\newenvironment{rem}{%
\bigskip
\noindent \textsl{{\sl Remark. }}}{\bigskip}
\newenvironment{rems}{%
\bigskip
\noindent \textsl{{\sl Remarks. }}}{\bigskip}

\newcounter {own}
\def\theown {\thesection.\arabic{own}}

\newenvironment{pf}[1][]{%
 \vskip 3mm
 \noindent
 \ifthenelse{\equal{#1}{}}%
  {{\slshape Proof. }}%
  {{\slshape #1.} }%
 }%
{\qed\bigskip}

\newcounter{alphabet}
\newcounter{tmp}

\begin{document}
\bibliographystyle{amsplain}
\newcommand{\co}{{\overline{\operatorname{co}}}}
\newcommand{\A}{{\mathcal A}}
\newcommand{\B}{{\mathcal B}}
\newcommand{\T}{{\mathcal T}}
\newcommand{\U}{{\mathcal U}}
\newcommand{\es}{{\mathcal S}}
\newcommand{\LU}{{\mathcal{LU}}}
\newcommand{\ZF}{{\mathcal{ZF}}}
\newcommand{\IR}{{\mathbb R}}
\newcommand{\IC}{{\mathbb C}}
\newcommand{\IN}{{\mathbb N}}
\newcommand{\K}{{\mathcal K}}
\newcommand{\X}{{\mathcal X}}
\newcommand{\PP}{{\mathcal P}}
\newcommand{\uhp}{{\mathbb H}}
\newcommand{\Z}{{\mathbb Z}}
\newcommand{\N}{{\mathcal N}}
\newcommand{\M}{{\mathcal M}}
\newcommand{\SCC}{{\mathcal{SCC}}}
\newcommand{\CC}{{\mathcal C}}
\newcommand{\st}{{\mathcal{SS}}}
\newcommand{\D}{{\mathbb D}}
\newcommand{\sphere}{{\widehat{\mathbb C}}}
\newcommand{\image}{{\operatorname{Im}\,}}
\newcommand{\Aut}{{\operatorname{Aut}}}
\newcommand{\real}{{\operatorname{Re}\,}}
\newcommand{\kernel}{{\operatorname{Ker}}}
\newcommand{\ord}{{\operatorname{ord}}}
\newcommand{\id}{{\operatorname{id}}}
\newcommand{\mob}{{\text{\rm M\"{o}b}}}
\newcommand{\Int}{{\operatorname{Int}\,}}
\newcommand{\Sign}{{\operatorname{Sign}}}
\newcommand{\diam}{{\operatorname{diam}}}
\newcommand{\inv}{^{-1}}
\newcommand{\area}{{\operatorname{Area}}}
\newcommand{\eit}{{e^{i\theta}}}
\newcommand{\dist}{{\operatorname{dist}}}
\newcommand{\arctanh}{{\operatorname{arctanh}}}
\newcommand{\remark}{\vskip .3cm \noindent {\sl Remark.} \@}
\newcommand{\remarks}{\vskip .3cm \noindent {\sl Remarks.} \@}
\newcommand{\ucv}{{\operatorname{UCV}}}
%%%%%%%%%%%%%%%%%%%%%%%%%%%%%%%%%%%%%%%%%%%%%%%%%%%%%%%%%%%%%%%%%%%%%%%%%%%%%%%%%5
%%%%%%%%%%%%%%%%%%%%%%%%%%%%%%%%%%%%%%%%%%%%%%%%%%%%%%%%%%%%%%%%%%%%%%
\def\be{\begin{equation}}
\def\ee{\end{equation}}
\newcommand{\sep}{\itemsep -0.01in}
\newcommand{\seps}{\itemsep -0.02in}
\newcommand{\sepss}{\itemsep -0.03in}
\newcommand{\bee}{\begin{enumerate}}
\newcommand{\eee}{\end{enumerate}}
\newcommand{\pays}{\!\!\!\!}
\newcommand{\pay}{\!\!\!}
\newcommand{\blem}{\begin{lem}}
\newcommand{\elem}{\end{lem}}
\newcommand{\bthm}{\begin{thm}}
\newcommand{\ethm}{\end{thm}}
\newcommand{\bcor}{\begin{cor}}
\newcommand{\ecor}{\end{cor}}
\newcommand{\beg}{\begin{example}}
\newcommand{\eeg}{\end{example}}
\newcommand{\begs}{\begin{examples}}
\newcommand{\eegs}{\end{examples}}
\newcommand{\bdefe}{\begin{defn}}
\newcommand{\edefe}{\end{defn}}
\newcommand{\bprob}{\begin{prob}}
\newcommand{\eprob}{\end{prob}}
\newcommand{\bcon}{\begin{conj}}
\newcommand{\econ}{\end{conj}}
\newcommand{\bprop}{\begin{prop}}
\newcommand{\eprop}{\end{prop}}
\newcommand{\br}{\begin{rem}}
\newcommand{\er}{\end{rem}}
\newcommand{\brs}{\begin{rems}}
\newcommand{\ers}{\end{rems}}
\newcommand{\bo}{\begin{obser}}
\newcommand{\eo}{\end{obser}}
\newcommand{\bos}{\begin{obsers}}
\newcommand{\eos}{\end{obsers}}
\newcommand{\bpf}{\begin{pf}}
\newcommand{\epf}{\end{pf}}
\newcommand{\ba}{\begin{array}}
\newcommand{\ea}{\end{array}}
\newcommand{\beq}{\begin{eqnarray}}
\newcommand{\beqq}{\begin{eqnarray*}}
\newcommand{\eeq}{\end{eqnarray}}
\newcommand{\eeqq}{\end{eqnarray*}}
\newcommand{\llra}{\longleftrightarrow}
\newcommand{\lra}{\longrightarrow}
\newcommand{\lla}{\longleftarrow}
\newcommand{\Llra}{\Longleftrightarrow}
\newcommand{\Lra}{\Longrightarrow}
\newcommand{\Lla}{\Longleftarrow}
\newcommand{\Ra}{\Rightarrow}
\newcommand{\La}{\Leftarrow}
\newcommand{\ra}{\rightarrow}
\newcommand{\la}{\leftarrow}

\title[Geometric properties of Clausen's  Hypergeometric Function $_3F_2(a,b,c;d,e;z)$]
{Geometric properties of Clausen's  Hypergeometric Function $_3F_2(a,b,c;d,e;z)$}

\author[K.Chandrasekran]{Koneri Chandrasekran}
\address{Koneri Chandrasekran \\ Department of Mathematics \\ Jeppiaar SRR Engineering College, Affiliated to Anna University \\ Chennai 603 103, India}
\email{kchandru2014@gmail.com}

\author[D. J. Prabhakaran ]{Devasir John Prabhakaran}
\address{Devasir John Prabhakaran \\ Department of Mathematics \\ MIT Campus, Anna University \\ Chennai 600 044, India}
\email{asirprabha@yahoo.com}


\begin{abstract}
The Clausen's Hypergeometric Function is given by
\beqq
_3F_2(a,b,c;d,e;z)=\sum_{n=0}^{\infty}\frac{(a)_n(b)_n(c)_n}{(d)_n(e)_n(1)_n}z^n\, ; \, \, \, a,b,c,d,e\in \IC
\eeqq provided $d,\, e\, \neq 0,-1,-2,\cdots$ in  unit disc $\D =\{z\in \IC\,:\, |z|<1\}$. In this paper, an operator $\mathcal{I}_{a,b,c}(f)(z)$ involving Clausen's Hypergeometric Function by means of Hadamard  Product is introduced. Geometric properties of $\mathcal{I}_{a,b,c}(f)(z)$ are obtained based on its Taylor's co-efficient.
\end{abstract}
\maketitle

\section{Introduction and preliminaries}
Let $\mathcal{A}$ be the class of functions
\beq\label{inteq0}
f(z)= z+\sum_{n=2}^{\infty}\, a_n\,z^n
\eeq analytic in the open unit disc $\D =\{z\in \IC\,:\, |z|<1\}$ of the complex plane.\\

Let   $\es$, \, $\es^{*}$, $\mathcal{C}$ and $\mathcal{K}$ be the class of univalent , Starlike, Convex and Close-to-Convex functions respectively. We will be particularly focussing on the classes $\es^{*}_{\lambda},\,  \lambda >0$ and $\mathcal{C}_{\lambda}$ are defined by $$\mathcal{S}^{*}_{\lambda}\, =\, \left\{f\in \mathcal{A}\, |\, \left|\frac{zf'(z)}{f(z)}-1\right| \, < \, \lambda,\, z\in \D \right\}.$$ and $$\mathcal{C_{\lambda}}=\left\{f\in \mathcal{A}\, |\, zf'(z)\in \mathcal{S}^{*}_{\lambda}\right\}.$$

The following  are the sufficient conditions for which the function $f$ is in $\mathcal{S}^{*}$  and $\mathcal{C}_{\lambda}$ respectively

\beq\label{inteq2}
\displaystyle \sum_{n=2}^{\infty}(n+\lambda-1)|a_n| \leq \lambda.
\eeq
and
\beq\label{inteq}
\displaystyle \sum_{n=2}^{\infty}\, n\, (n+\lambda-1)|a_n| \leq \lambda.
\eeq

Let for $\beta <1$, $${\mathcal R}(\beta)=\{f\in {\mathcal A}:  \exists \  \eta \in \left(-\frac{\pi}{2}, \frac{\pi}{2}\right)\, |\, {\rm Re}\, \ [ e^{i\eta}(f'(z)-\beta)] > 0, \quad z\in\D \}. $$

Note that when $\beta \ge 0$, we have ${\mathcal R}(\beta) \subset {\mathcal S}$ and for each $\beta <0 ,\ \  {\mathcal R}(\beta) $ contains also non univalent functions.\\

The concept of uniformly convex and uniformly starlike functions were introduced by Goodman\cite{Good-1991-Ann-PM,Good-1991-JMAA} and denoted by UCV and UST respectively. Subsequently R{\o}nning \cite{Ronn-1993-Proc-ams} and Ma and Minda \cite{Ma-Minda-1992-Ann-PM} independently gave the one variable analytic characterization of the class UCV as follows: $f\in UCV$ if and only if $\displaystyle Re\left(1+\frac{zf^{''}(z)}{f'(z)}\right) > \left|\frac{zf''(z)}{f'(z)}\right|$.\\

In \cite{Good-1991-Ann-PM}, the condition on the function  $f$ defined in (\ref{inteq0}) belongs to $UCV$ is that
\beq\label{lem4eq1}
\sum_{n=2}^{\infty}\,n\, (n-1)|a_n|\leq \frac{1}{3}
\eeq

The  subclass $\es_p$ of starlike functions was introduced by R{\o}nning \cite{Ronn-1993-Proc-ams} in the following way
\beq\label{inteq3}
\displaystyle \es_p = \{ F \in \es^{*}/F(z)=zf'(z),\, f(z) \in UCV\}.
\eeq

It is easily seen that $f(z)\in \es_p$ if and only if
\beq\label{inteq4}
\left|\frac{zf'(a)}{f(z)}-1\right| < Re\left(\frac{zf'(z)}{f(z)}\right),\, z\in \D
\eeq
 and has a sufficient condition
 \beq \label{lem2eq1}
\sum_{n=2}^{\infty}(2n-1)|a_n|\leq 1,
\eeq

We observe that $\es_P$ is the class of functions for which the domain of values of $\displaystyle zf'(z)/f(z),$ $ z \in \D$ is the region $\Omega$ defined by $\Omega =\{ w : Re(w) > |w-1| \}$. Note that $\Omega$ is the interior of a parabola in the right half plane which is symmetric about the real axis and has vertex at $\left(\frac{1}{2}, 0 \right)$. It is well- known that the function $$\phi(z)=1+\frac{2}{\pi^2}\left(log \frac{1+\sqrt{z}}{1-\sqrt{z}}\right)^2$$ maps unit disc $\D$ onto the parabolic region $\Omega$ and hence is in $\es_p.$\\

Let $ \displaystyle f(z)= z+\sum_{n=2}^{\infty}\, a_n\,z^n $ and  $ \displaystyle g(z)= z+\sum_{n=2}^{\infty}\, b_n\,z^n $ be analytic in $\D$. Then the Hadamard product or convolution of $f(z)$  and $g(z)$ is defined by $$f(z)*g(z)= z+\sum_{n=2}^{\infty} a_nb_n z^n.$$

For any complex variable $a$, define the ascending factorial notation $$(a)_n = a(a+1)(a+2)\cdots(a+n-1) = a (a+1)_{n-1}$$ for $n\geq1$ and $(a)_0=1$ for $a\neq0$. When $a$ is neither zero nor a negative integer, we have $(a)_n = \Gamma(n+a)/\Gamma(a).$\\

The Clausen's hypergeometric function $_3F_2(a,b,c;d,e;z)$ is defined by
\beq\label{inteq5}
_3F_2(a,b,c;d,e;z)=\sum_{n=0}^{\infty}\frac{(a)_n(b)_n(c)_n}{(d)_n(e)_n(1)_n}z^n;\, \, \, a,b,c,d,e\in \IC
\eeq provided $d,\, e\, \neq 0,-1,-2,-3\cdots$ Which is an analytic function in unit disc $\D$.\\

For $n=1$ in Theorem 1 in Miller and Paris \cite{Mill-Paris-2012-ITSF} yields the following formula
\beq\label{inteq6}
_3F_2(a,b,c;b+1,c+1;1)&=& \frac{bc}{c-b}\Gamma(1-a)\left[\frac{\Gamma(b)}{\Gamma(1-a+b)}-\frac{\Gamma(c)}{\Gamma(1-a+c)}\right],
\eeq provided that $Re(2- a) > 0$,  $b > a-1$ and $c > a-1$. Alternatively we drive the same by putting $m=n=0$ in the equation (3) in Shpot and Srivastava \cite{Shpot-Srivas-2015-AMC}. We use the formula (\ref{inteq6}) to prove our main results. \\

For $f\in \mathcal{A}$, we define the operator $\mathcal{I}_{a,b,c}(f)(z)$
\beq\label{inteq7}
\mathcal{I}_{a,b,c}(f)(z) = z\, _3F_2(a,b,c;b+1,c+1;z)*f(z) = z+\sum_{n=2}^{\infty} A_n\, z^n
\eeq
with $A_1=1$ and for $n\geq 1,$
\beqq
A_n&=&\frac{(a)_{n-1}(b)_{n-1}(c)_{n-1}}{(b+1)_{n-1}(c+1)_{n-1}(1)_{n-1}}\, a_n.\\
\eeqq

The following Lemma is useful to prove our main results.
\blem \label{lem1eq1}
Let $a,b,c > 0$. Then we have the following
\begin{enumerate}
\item For $ b, c > a-1$
\begin{eqnarray*}
\sum_{n=0}^{\infty} \frac{(n+1)(a)_n\, (b)_n\, (c)_n }
{(b+1)_n\, (c+1)_n\, (1)_n} &=&  \frac{bc\, \Gamma(1-a) }{c-b}\left[\frac{(1-b)\Gamma(b)}{\Gamma(1-a+b)}-\frac{(1-c)\Gamma(c)}{\Gamma(1-a+c)}\right].
\end{eqnarray*}
\item For $ b, c > a-1$
\begin{eqnarray*}
\sum_{n=0}^{\infty} \frac{(n+1)^2(a)_n\, (b)_n\, (c)_n }
{(b+1)_n\, (c+1)_n\, (1)_n} &=&  \frac{bc\, \Gamma(1-a)}{c-b}\left[\frac{(1-b)^2\Gamma(b)}{\Gamma(1-a+b)}-\frac{(1-c)^2\Gamma(c)}{\Gamma(1-a+c)}\right].
\end{eqnarray*}
\item For $ b, c > a-1$. Then
\begin{eqnarray*}
\sum_{n=0}^{\infty} \frac{(n+1)^3(a)_n\, (b)_n\, (c)_n }
{(b+1)_n\, (c+1)_n\, (1)_n} &=&  \frac{bc\, \Gamma(1-a)}{c-b} \left[\frac{[(1-b)^3-b^2]\Gamma(b)}{\Gamma(1-a+b)}-\frac{[(1-c)^3-c^2]\Gamma(c)}{\Gamma(1-a+c)}\right].
\end{eqnarray*}
\item For $a\neq 1,\, b\neq 1,\,$ and $c\neq 1$ with$\, b, c > Max\{0,  a-1\}$
\begin{eqnarray*}
\sum_{n=0}^{\infty} \frac{(a)_n\, (b)_n\, (c)_n }
{(b+1)_n\, (c+1)_n\, (1)_{(n+1)}} &=&  \frac{bc}{(a-1)(b-1)(c-1)}\\ && \times\left[\frac{\Gamma(2-a)}{c-b}\left(\frac{(c-1)\Gamma(b)}{\Gamma(1-a+b)}-\frac{(b-1)\Gamma(c)}{\Gamma(1-a+c)}\right)-1\right].
\end{eqnarray*}
\end{enumerate}
\elem
\bpf (1) Using ascending factorial notation, we can write
\begin{eqnarray*}
\sum_{n=0}^{\infty} \frac{(n+1)(a)_n\, (b)_n\, (c)_n }
{(b+1)_n\, (c+1)_n\, (1)_n}
&=& \sum_{n=0}^{\infty} \frac{(a)_{n+1}\, (b)_{n+1}\, (c)_{n+1} }
{(b+1)_{n+1}\, (c+1)_{n+1}\, (1)_{n-1}} +\sum_{n=0}^{\infty} \frac{(a)_n\, (b)_n\, (c)_n }
{(b+1)_n\, (c+1)_n\, (1)_n}
\end{eqnarray*}
Using the formula (\ref{inteq6}) and the fact that $\Gamma(1-a)=-a\Gamma(-a)$, the above reduces to
\begin{eqnarray*}
\sum_{n=0}^{\infty} \frac{(n+1)(a)_n\, (b)_n\, (c)_n }
{(b+1)_n\, (c+1)_n\, (1)_n} &=& \frac{bc}{c -b}\, \Gamma(1-a)\,\left[\frac{(1-b)\Gamma(b)}{\Gamma(1-a+b)}-\frac{(1-c)\Gamma(c)}{\Gamma(1-a+c)}\right]
\end{eqnarray*}
Hence, (1) is proved.\\

(2) Using ascending factorial notation and $ (n+1)^2 = n(n-1)+3n+1$, we can write\\
\begin{flushleft}
$ \displaystyle \sum_{n=0}^{\infty} \frac{(n+1)^2(a)_n\, (b)_n\, (c)_n }{(b+1)_n\, (c+1)_n\, (1)_n} $
\end{flushleft}
\begin{eqnarray*}
&=& \sum_{n=2}^{\infty} \frac{(a)_{n}\, (b)_{n}\, (c)_{n} }{(b+1)_{n}\, (c+1)_{n}\, (1)_{n-2}} +\sum_{n=1}^{\infty} \frac{3 \, (a)_n\, (b)_n\, (c)_n }
{(b+1)_n\, (c+1)_n\, (1)_{n-1}}  +\sum_{n=0}^{\infty} \frac{(a)_n\, (b)_n\, (c)_n }{(b+1)_n\, (c+1)_n\, (1)_{n}}
\end{eqnarray*}
Using the formula (\ref{inteq6})
and  $\Gamma(1-a)=-a\Gamma(-a)$, the above reduces to
\begin{eqnarray*}
\sum_{n=0}^{\infty} \frac{(n+1)^2(a)_n\, (b)_n\, (c)_n }{(b+1)_n\, (c+1)_n\, (1)_n}&=& \frac{bc\, \Gamma(1-a)}{c-b}\times  \left[\frac{(1-b)^2\Gamma(b)}{\Gamma(1-a+b)}-\frac{(1-c)^2\Gamma(c)}{\Gamma(1-a+c)}\right]
\end{eqnarray*}
Which completes the proof of (2).\\

(3)   Using shifted factorial notation and by adjusting coefficients suitably, we can write\\
\begin{flushleft}
$ \displaystyle \sum_{n=0}^{\infty} \frac{(n+1)^3(a)_n\, (b)_n\, (c)_n }{(b+1)_n\, (c+1)_n\, (1)_n} $
\end{flushleft}
\begin{eqnarray*}
&=&\sum_{n=0}^{\infty} \frac{(a)_{n+3}\, (b)_{n+3}\, (c)_{n+3} }
{(b+1)_{n+3}\, (c+1)_{n+3}\, (1)_{n}}  +5\sum_{n=0}^{\infty} \frac{(a)_{n+2}\, (b)_{n+2}\, (c)_{n+2} }
{(b+1)_{n+2}\, (c+1)_{n+2}\, (1)_{n}}  \\&& +6\sum_{n=0}^{\infty} \frac{(a)_{n+1}\, (b)_{n+1}\, (c)_{n+1}}
{(b+1)_{n+1}\, (c+1)_{n+1}\, (1)_{n}} +\sum_{n=0}^{\infty} \frac{(a)_n\, (b)_n\, (c)_n }{(b+1)_n\, (c+1)_n\, (1)_{n}}
\end{eqnarray*}
Using the formula (\ref{inteq6}),  we have
\begin{eqnarray*}
\sum_{n=0}^{\infty} \frac{(n+1)^3(a)_n\, (b)_n\, (c)_n }{(b+1)_n\, (c+1)_n\, (1)_n}&=& \frac{bc\, \Gamma(1-a)}{c-b}\times  \left[\frac{[(1-b)^3-b^2]\Gamma(b)}{\Gamma(1-a+b)}-\frac{[(1-c)^3-c^2]\Gamma(c)}{\Gamma(1-a+c)}\right]
\end{eqnarray*}
and the conclusion follows.\\

(4)  Let $a$ be a positive real number such that $a\neq 1$, $b\neq 1$ and $c\neq1$ with $b, c > Max\{0, a-1\}$. We find that \\

$\displaystyle\sum_{n=0}^{\infty} \frac{(a)_n\, (b)_n\, (c)_n }{(b+1)_n\, (c+1)_n\, (1)_{n+1}} $
\begin{eqnarray*}
&=& \frac{bc}{(a-1)(b-1)(c-1)} \left[\frac{\Gamma(2-a)}{c-b}\left(\frac{(c-1)\Gamma(b)}{\Gamma(1-a+b)}-\frac{(b-1)\Gamma(c)}{\Gamma(1-a+c)}\right)-1\right]
\end{eqnarray*}
and the result follows.
\epf

\section{Starlikeness of $z _3F_2(a,b,c;b+1,c+1;z)$}
\bthm\label{thm1eq0}
 Let $a \in {\Bbb C} \backslash \{ 0 \} $,\, $b > |a|-1 $   and $c > |a|-1 .$ The sufficient condition for the function $z _3F_2(a,b,c;b+1,c+1;z) $ belongs to the class $ \es^{*}_{\lambda}, \,  0 < \lambda  \leq 1 $ is that
\beq\label{thm1eq1}
   \frac{bc}{c-b}\, \Gamma(1-|a|)\left[\frac{\Gamma(b)\, (\lambda-b)}{\Gamma(1-|a|+b)}-\frac{\Gamma(c)\, (\lambda-c)}{\Gamma(1-|a|+c)}\right] &\leq& 2\lambda
\eeq
\ethm
\bpf  Let $f(z)=z _3F_2(a,b,c;b+1,c+1;z)$. Then by the equation (\ref{inteq2}), it is enough to show that
\begin{eqnarray*}
  T &=& \sum_{n=2}^{\infty}(n+\lambda-1)|a_n|\leq \lambda.
\end{eqnarray*}
  Since $ f \in \es$, we have $|a_n| \leq 1$, and using the fact that  $|(a)_n|\leq  (|a|)_n$,
\begin{eqnarray*}
  T &\leq&  \sum_{n=2}^{\infty} (n-1+\lambda)\left(\frac{(|a|)_{n-1}(b)_{n-1}(c)_{n-1}}{(b+1)_{n-1}(c+1)_{n-1}(1)_{n-1}}\right)
\end{eqnarray*}
Using the formula (\ref{inteq6}) and the result (1) of Lemma \ref{lem1eq1}   in above equation, we get
\begin{eqnarray*}
T &\leq& \frac{bc}{c-b}\, \Gamma(1-|a|)\left[\frac{\Gamma(b)(\lambda-b)}{\Gamma(1-|a|+b)}-\frac{\Gamma(c)(\lambda-c)}{\Gamma(1-|a|+c)}\right]  -\lambda
\end{eqnarray*}
Because of (\ref{thm1eq1}), the above expression is bounded by $\lambda$ and hence
\begin{eqnarray*}
  T &\leq &  \frac{bc}{c-b}\, \Gamma(1-|a|)\left[\frac{\Gamma(b)(\lambda-b)}{\Gamma(1-|a|+b)}-\frac{\Gamma(c)(\lambda-c)}{\Gamma(1-|a|+c)}\right]  -\lambda \leq \lambda
\end{eqnarray*}
Hence $z _3F_2(a,b,c;b+1,c+1;z) $ belongs to the class $ \es^{*}_{\lambda}. $
\epf
\bthm\label{thm2eq001}
 Let $a \in {\Bbb C} \backslash \{ 0 \} ,\, c > 0, \, b > 0, \, |a|\neq1,\, b \neq1,\, c \neq 1,\,  b > |a|-1 $   and $c > |a|-1 .$ For  $ 0 < \lambda \leq 1$, assume that
 \beq\label{thm2eq1}
   \frac{bc}{c-b}\, \Gamma(1-|a|)\left[\left(\frac{b-\lambda}{b-1}\right)\frac{\Gamma(b)}{\Gamma(1-|a|+b)}-\left(\frac{c-\lambda}{c-1}\right)\frac{\Gamma(c)}{\Gamma(1-|a|+c)}\right] \\ \leq \lambda\left(1+\frac{1}{2(1-\beta)}\right)+\frac{(\lambda-1)bc}{(|a|-1)(b-1)(c-1)}\nonumber
\eeq
 Then the integral  operator $\mathcal{I}_{a,b,c}(f)$ maps $ \mathcal{R}(\beta)$ into $\es^{*}_{\lambda}$.
\ethm
\bpf
 Let $a \in {\Bbb C} \backslash \{ 0 \} ,\, c > 0, \, b > 0, \, |a|\neq1,\, b \neq1,\, c \neq 1,\,  b > |a|-1 $   and $c > |a|-1 .$\\

Suppose that $\displaystyle f(z) $ is defined in (\ref{inteq0}) is in $ \mathcal{R}(\beta)$. By MacGregor \cite{MacGregor-1962-Trans-ams}, We have
\beq\label{thm2eq2}
|a_n|\leq \frac{2(1-\beta)}{n}.
\eeq
Consider the integral operator $\mathcal{I}_{a,b,c}(f)$ is defined by (\ref{inteq7}). According to the equation (\ref{inteq2}), we need to show that
\begin{eqnarray*}
  T &=& \sum_{n=2}^{\infty}(n+\lambda-1)|A_n|\leq \lambda.
\end{eqnarray*}
Then, we have
\begin{eqnarray*}
  T &=&  \sum_{n=2}^{\infty} [(n-1)+\lambda] \left|\frac{(a)_{n-1}(b)_{n-1}(c)_{n-1}}{(b+1)_{n-1}(c+1)_{n-1}(1)_{n-1}}\right| |a_n|\\
\end{eqnarray*}
Using (\ref{thm2eq2})  in above, we have
\begin{eqnarray*}
  T   &\leq&  2(1-\beta)\left[\sum_{n=1}^{\infty} (n+1) \frac{(|a|)_{n}(b)_{n}(c)_{n}}{(b+1)_{n}(c+1)_{n}(1)_{n}(n+1)} \right.\\ && \left.+ (\lambda-1)\sum_{n=1}^{\infty}  \frac{(|a|)_{n}(b)_{n}(c)_{n}}{(b+1)_{n}(c+1)_{n}(1)_{n}}\left(\frac{1}{n+1}\right)  \right]:=T_1
\end{eqnarray*}
Using the formula given by (\ref{inteq6}) and the results (1) and (4) of Lemma \ref{lem1eq1}, we find that
\begin{eqnarray*}
  T_1 &\leq&  2(1-\beta)\left[ \frac{bc}{c-b}\Gamma(1-|a|)\left(\frac{\Gamma(b)}{\Gamma(1-|a|+b)}-\frac{\Gamma(c)}{\Gamma(1-|a|+c)}\right)-1  \right.\\ && \left.+ (\lambda-1)\left(\frac{bc}{(a-1)(b-1)(c-1)}\left(\left(\frac{2-a}{c-b}\right)\left(\frac{(c-1)\Gamma(b)}{\Gamma(1-|a|+b)}-\frac{(b-1)\Gamma(c)}{\Gamma(1-|a|+c)}\right) \right.\right. \right.\\ && \left.\left. \left.-1 \right)-1 \right) \right]
\end{eqnarray*}
Using the fact that $\Gamma(a+1)=a\Gamma(a)$, we have
\begin{eqnarray*}
T_1 &\leq&  2(1-\beta)\left[ \frac{bc}{c-b}\Gamma(1-|a|)\left(\frac{\Gamma(b)}{\Gamma(1-|a|+b)}\left(\frac{b-\lambda}{b-1}\right)-\frac{\Gamma(c)}{\Gamma(1-|a|+c)} \left(\frac{c-\lambda}{c-1}\right)\right)\right.\\ && \left.
+ \frac{(\lambda-1)\, bc}{(a-1)(b-1)(c-1)}-\lambda\right]
\end{eqnarray*}
Under the condition given by  (\ref{thm2eq1})
\begin{eqnarray*}
 2(1-\beta)\left[ \frac{bc}{c-b}\Gamma(1-|a|)\left(\frac{\Gamma(b)}{\Gamma(1-|a|+b)}\left(\frac{b-\lambda}{b-1}\right)-\frac{\Gamma(c)}{\Gamma(1-|a|+c)} \left(\frac{c-\lambda}{c-1}\right)\right)\right.&&\\  \left.
+ \frac{(\lambda-1)\, bc}{(a-1)(b-1)(c-1)}-\lambda\right] &\leq& \lambda
\end{eqnarray*}
Thus, the inequality $T \leq T_1 \leq \lambda $ and (\ref{thm2eq2}) are holds. Therefore, we conclude that the operator $\mathcal{I}_{a,b,c}(f)$ maps $ \mathcal{R}(\beta)$ into $\es^{*}_{\lambda}$. Which completes the proof of the theorem.
\epf

When $\lambda =1 $, we get the following results from Theorem \ref{thm2eq001}.
\bcor
Let $a \in {\Bbb C} \backslash \{ 0 \},\,  b > |a|-1 $   and $c > |a|-1 .$ Assume that
 \beq\label{cor2eq1}
   \frac{bc}{c-b}\, \Gamma(1-|a|)\left[\frac{\Gamma(b)}{\Gamma(1-|a|+b)}-\frac{\Gamma(c)}{\Gamma(1-|a|+c)}\right]  \leq 1+\frac{1}{2(1-\beta)}\nonumber
\eeq
 Then the integral  operator $\mathcal{I}_{a,b,c}(f)$ maps $ \mathcal{R}(\beta)$ into $\es^{*}_{1}$.
\ecor
\bthm\label{thm3eq0}  Let $a \in {\Bbb C} \backslash \{ 0 \} ,\,  b > |a|-1$ and $c > |a|-1$. Suppose that $a$ and $0 < \lambda \leq 1$ satisfy the condition
 \beq\label{thm3eq1}
 \frac{bc}{c-b}\, \Gamma(1-|a|)\left(\frac{(1-b)(\lambda-b)\Gamma(b)}{\Gamma(1-|a|+b)}-\frac{(1-c)(\lambda-c)\Gamma(c)}{\Gamma(1-|a|+c)}\right)  &\leq& 2\lambda
\eeq
then the integral  operator $\mathcal{I}_{a,b,c}(f)$ maps $\mathcal{S}$ to $\es^{*}_{\lambda}$.
\ethm
\bpf Let $a \in {\Bbb C} \backslash \{ 0 \} , \, b > |a|-1$ and $c > |a|-1$. Suppose that $\displaystyle f(z)$ is defined by (\ref{inteq0}) is in $\es$, then we know that
\beq\label{thm3eq2}
|a_n|\leq n.
\eeq
Suppose that the integral  operator $\mathcal{I}_{a,b,c}(f)$ is defined by (\ref{inteq7}).
In the view of the equation (\ref{inteq2}), it is enough to show that
\begin{eqnarray*}
  T &=& \sum_{n=2}^{\infty}(n+\lambda-1)|A_n|\leq \lambda.
\end{eqnarray*}
Using the fact that $|(a)_n|\leq  (|a|)_n$ and the equation ($\ref{thm3eq2}$) in above, we have
\begin{eqnarray*}
  T &= & \sum_{n=2}^{\infty} (n-1+\lambda) \left|\frac{(a)_{n-1}(b)_{n-1}(c)_{n-1}}{(b+1)_{n-1}(c+1)_{n-1}(1)_{n-1}}\right||a_n|
\end{eqnarray*}
Using (1) and (2) of Lemma \ref{lem1eq1}, we find that
\begin{eqnarray*}
T &\leq&  \frac{bc}{c-b}\Gamma(1-|a|)\left(\frac{(1-b)\Gamma(b)}{\Gamma(1-|a|+b)}[\lambda-b]-\frac{(1-c)\Gamma(c)}{\Gamma(1-|a|+c)}[\lambda -c]\right)  -\lambda
\end{eqnarray*}
By the condition given by (\ref{thm3eq1}), the above expression is bounded by $\lambda$ and hence
\begin{eqnarray*}
 \frac{bc}{c-b}\Gamma(1-|a|)\left(\frac{(1-b)(\lambda-b)\Gamma(b)}{\Gamma(1-|a|+b)}-\frac{(1-c)(\lambda-c)\Gamma(c)}{\Gamma(1-|a|+c)}\right) -\lambda &\leq& \lambda
\end{eqnarray*}
Under the stated condition, The integral operator $\mathcal{I}_{a,b,c}(f)$ maps $\es$ into $\es^{*}_{\lambda}$. Which gives the appropriate conclusion.
\epf
\section{Convexity of $z _3F_2(a,b,c;b+1,c+1;z)$ }
\bthm\label{thm10eq1}
 Let $a \in {\Bbb C} \backslash \{ 0 \} $,\, $b > |a|-1 $, $c > |a|-1 $ and $ 0 < \lambda \leq 1$. The sufficient condition for the function $z _3F_2(a,b,c;b+1,c+1;z) $ belongs to the class $ \mathcal{C}_{\lambda}$ is that
 \beq\label{thm10eq10}
    \frac{bc}{c-b}\, \Gamma(1-|a|)\left[\frac{(1-b)\,(\lambda-b)\, \Gamma(b)}{\Gamma(1-|a|+b)}-\frac{(1-b)(\lambda-c)\Gamma(c)}{\Gamma(1-|a|+c)}\right]\leq 2\lambda.\nonumber
   \eeq
\ethm
\bpf The proof is similar to Theorem \ref{thm3eq0}. So we omit the details.
\epf

\bthm\label{thm11eq0}  Let $a \in {\Bbb C} \backslash \{ 0 \} ,\, b > 0,\, c > 0,\,   |a| \neq 1 ,\, b\neq 1$,\, $c \neq 1, \, b > |a|-1$,\,  $c > |a|-1$ and $ 0 < \lambda \leq 1$.  For $0 \leq \beta <1 $, assume that
 \beq\label{thm11eq1}
 \frac{bc}{c-b}\Gamma(1-|a|)\left(\frac{(\lambda-b)\, \Gamma(b)}{\Gamma(1-|a|+b)}-\frac{(\lambda-c)\,\Gamma(c)}{\Gamma(1-|a|+c)}\right)  &\leq& \lambda\left( \frac{1}{2(1-\beta)}+1\right)\nonumber
\eeq
then the operator $\mathcal{I}_{a,b,c}(f)$ maps $\mathcal{R}(\beta)$ into $ \mathcal{C}_{\lambda} $.
\ethm
\bpf The proof is similar to Theorem \ref{thm1eq0}. So we omit the details.
\epf
\bthm\label{thm12eq0}  Let $a \in {\Bbb C} \backslash \{ 0 \},\,  b > |a|-1$ and $c > |a|-1$. Suppose that $a$ and $ 0 < \lambda \leq 1$ satisfy the condition
 \beq\label{thm12eq1}
 \frac{bc}{c-b}\Gamma(1-|a|)\left(\frac{[b^2-b-b^3+\lambda(1-b)^2]\Gamma(b)}{\Gamma(1-|a|+b)}\right. \\  \left.-\frac{[c^2-c-c^3+\lambda(1-c)^2]\Gamma(c)}{\Gamma(1-|a|+c)}\right) &\leq& 2\lambda\nonumber
\eeq
then $\mathcal{I}_{a,b,c}(f)$ maps $\es$ into $ \mathcal{C}_{\lambda} $.
\ethm
\bpf Let $a \in {\Bbb C} \backslash \{ 0 \} , \, b > |a|-1$ and $c > |a|-1$. Suppose that $\displaystyle f(z)$ is defined by (\ref{inteq0}) is in $\es$, then
\beq\label{thm12eq2}
|a_n|\leq n.
\eeq
Suppose the integral operator $\mathcal{I}_{a,b,c}(f)$ defined by (\ref{inteq7}). In the view of the sufficient condition given by (\ref{inteq}), it is enough to prove that
\begin{eqnarray*}
  T &=& \sum_{n=2}^{\infty}\,n\, (n+\lambda-1)\, |B_n|\leq \lambda .
\end{eqnarray*}
i.e.,
\begin{eqnarray*}
  T &=& \sum_{n=2}^{\infty}n\, (n+\lambda-1)\, \left|\frac{(a)_{n-1}(b)_{n-1}(c)_{n-1}}{(b+1)_{n-1}(c+1)_{n-1}(1)_{n-1}}\right|\, |a_n|\leq \lambda.
\end{eqnarray*}
Using the fact that $|(a)_n|\leq  (|a|)_n$ and the equation $(\ref{thm12eq2})$ in above, we have
\begin{eqnarray*}
  &\leq&  \sum_{n=0}^{\infty}  \frac{(n+1)^3\, (|a|)_{n}(b)_{n}(c)_{n}}{(b+1)_{n}(c+1)_{n}(1)_{n}}+(\lambda-1)\sum_{n=0}^{\infty}  \frac{(n+1)^2\, (|a|)_{n}(b)_{n}(c)_{n}}{(b+1)_{n}(c+1)_{n}(1)_{n}}-\lambda
\end{eqnarray*}
Using (2) and (3) of Lemma \ref{lem1eq1}, we find that
\begin{eqnarray*}
T &\leq&   \frac{2\, bc}{c-b}\Gamma(1-|a|)\left(\frac{[(1-b)^3-b^2]\Gamma(b)}{\Gamma(1-|a|+b)}-\frac{[(1-c)^3-c^2]\Gamma(c)}{\Gamma(1-|a|+c)}\right)  \\ &&+ \frac{bc (\lambda-1)}{c-b}\Gamma(1-|a|)\left(\frac{(1-b)^2\Gamma(b)}{\Gamma(1-|a|+b)}-\frac{(1-c)^2\Gamma(c)}{\Gamma(1-|a|+c)}\right) -\lambda
\end{eqnarray*}
By the equation (\ref{thm12eq1}), the above expression is bounded by $\lambda$ and hence
\begin{eqnarray*}
\frac{bc}{c-b}\Gamma(1-|a|)\left(\frac{[b^2-b-b^3+\lambda(1-b)^2]\Gamma(b)}{\Gamma(1-|a|+b)}-\frac{[c^2-c-c^3+\lambda(1-c)^2]\Gamma(c)}{\Gamma(1-|a|+c)}\right) &\leq& 2\lambda
\end{eqnarray*}
Hence, the integral operator $\mathcal{I}_{a,b,c}(f)$ maps $\es$ into $\mathcal{C}_{\lambda}$ and the proof is complete.
\epf
\section{Admissibility condition of $z _3F_2(a,b,c;b+1,c+1;z)$ in Parabolic domain.}
\bthm\label{thm7eq1}
 Let $a \in {\Bbb C} \backslash \{ 0 \} $,\, $b > |a|-1 $   and $c > |a|-1 .$  The sufficient condition for the function $z _3F_2(a,b,c;b+1,c+1;z) $ belongs to the class $ UCV$ is that
 \beq\label{thm7eq10}
   \frac{bc}{c-b}\, \Gamma(1-|a|)\left[\frac{(1-c)\, \Gamma(c+1)\, }{\Gamma(1-|a|+c)}-\frac{(1-b)\,\Gamma(b+1)\, }{\Gamma(1-|a|+b)}\right] &\leq&  \frac{1}{3}.\nonumber
   \eeq
\ethm
\bpf
The proof is similar to Theorem \ref{thm3eq0}. So we omit the details.
\epf
\bthm\label{thm8eq0}  Let $a \in {\Bbb C} \backslash \{ 0 \} ,\, b > 0,\, c > 0,\,  b > |a|-1$ and $c > |a|-1$.  For $0 \leq \beta <1 $,  assume that
 \beq\label{thm8eq1}
  \frac{bc}{c-b}\Gamma(1-|a|)\left(\frac{\Gamma(c+1)}{\Gamma(1-|a|+c)}-\frac{\Gamma(b+1)}{\Gamma(1-|a|+b)}\right)  &\leq& \frac{1}{6(1-\beta)}
\eeq
then $I_{a,b,c}(f)$ maps $\mathcal{R}(\beta)$ into $ UCV $.
\ethm
\bpf Let $a \in {\Bbb C} \backslash \{ 0 \} ,\, b > 0,\, c > 0,\, b > |a|-1$ and $c > |a|-1$.\\

Consider the integral operator $\mathcal{I}_{a,b,c}(f)$  is given by (\ref{inteq7}). According to sufficient condition given by (\ref{lem4eq1}), it is enough to show that
\begin{eqnarray*}
  T &=& \sum_{n=2}^{\infty}n\, (n-1)\,|B_n|\leq \frac{1}{3}.
\end{eqnarray*}
Using the fact that $|(a)_n|\leq  (|a|)_n$ and the equation $(\ref{thm2eq2})$ in above, we have
\begin{eqnarray*}
  T &\leq&  2(1-\beta)\sum_{n=2}^{\infty} n\,(n-1)\, \left|\frac{(|a|)_{n-1}(b)_{n-1}(c)_{n-1}}{(b+1)_{n-1}(c+1)_{n-1}(1)_{n-1}n}\right|
 \end{eqnarray*}
Using the formula given by (\ref{inteq6}) and (1) of Lemma \ref{lem1eq1}, we find that
\begin{eqnarray*}
T &\leq&  2(1-\beta)\left[ \frac{bc}{c-b}\Gamma(1-|a|)\left(\frac{(1-b-1)\, \Gamma(b)}{\Gamma(1-|a|+b)}-\frac{(1-c-1)\,\Gamma(c)}{\Gamma(1-|a|+c)}\right) \right]
\end{eqnarray*}
Using the fact that $\Gamma(a+1)=a\Gamma(a)$, the above reduces to
\begin{eqnarray*}
T &\leq&  2(1-\beta)\left[ \frac{bc}{c-b}\Gamma(1-|a|)\left(\frac{\Gamma(c+1)}{\Gamma(1-|a|+c)}-\frac{\Gamma(b+1)}{\Gamma(1-|a|+b)}\right) \right]
\end{eqnarray*}
By the formula given by (\ref{thm8eq1}), the above expression is bounded by $\frac{1}{3}$ and hence
\begin{eqnarray*}
2(1-\beta)\left[ \frac{bc}{c-b}\Gamma(1-|a|)\left(\frac{\Gamma(c+1)}{\Gamma(1-|a|+c)}-\frac{\Gamma(b+1)}{\Gamma(1-|a|+b)}\right) \right] &\leq& \frac{1}{3}
\end{eqnarray*}
By Hypothesis, the operator $\mathcal{I}_{a,b,c}(f)$ maps $\mathcal{R}(\beta)$ into $UCV$ and the result follows.
\epf
\bthm\label{thm9eq0}  Let $a \in {\Bbb C} \backslash \{ 0 \} ,\,  b > |a|-1$ and $c > |a|-1$. Assume that
 \beq\label{thm9eq1}
  \frac{bc}{c-b}\Gamma(1-|a|)\left(\frac{(b^2-b-b^3)\Gamma(b)}{\Gamma(1-|a|+b)}-\frac{(c^2-c-c^3)\Gamma(c)}{\Gamma(1-|a|+c)}\right)  &\leq& \frac{1}{3}\nonumber
\eeq
then the integral  operator $\mathcal{I}_{a,b,c}(f)$ maps $\mathcal{S}$ to $UCV$.
\ethm
\bpf
The proof is similar to Theorem \ref{thm12eq0}. So we omit the details.
\epf
\section{Inclusion Properties of $z\, _3F_2(a,b,c;b+1,c+1;z) $  in $\es_p$-CLASS}
\bthm\label{thm4eq0}
 Let $a \in {\Bbb C} \backslash \{ 0 \} $,\, $b > |a|-1 $   and $c > |a|-1 .$  The sufficient condition for the function $z\, _3F_2(a,b,c;b+1,c+1;z) $ belongs to the class $ S_{p} $ is that
 \beq\label{thm4eq1}
   \frac{bc}{c-b}\, \Gamma(1-|a|)\left[\frac{(1-2b)\, \Gamma(b)}{\Gamma(1-|a|+b)}-\frac{(1-2c)\,\Gamma(c)}{\Gamma(1-|a|+c)}\right] &\leq& 2.\nonumber
\eeq
\ethm
\bpf
The proof is similar to Theorem \ref{thm8eq0}. So we omit the details.
\epf
\bthm\label{thm5eq0}  Let $a \in {\Bbb C} \backslash \{ 0 \} ,\, b > 0,\, c > 0,\,   |a| \neq 1 ,\, b\neq 1$,\, $c \neq 1, \, b > |a|-1$ and $c > |a|-1$. Assume that
 \beq\label{thm5eq1}
     \frac{bc}{c-b}\, \Gamma(1-|a|)\left(\left(\frac{2b-1}{b-1}\right)\frac{\Gamma(b)}{\Gamma(1-|a|+b)}-\left(\frac{2c-1}{c-1}\right)\frac{\Gamma(c)}{\Gamma(1-|a|+c)}\right)\\ +\frac{bc}{(|a|-1)(b-1)(c-1)} \leq \frac{1}{2(1-\beta)}+1.\nonumber
\eeq
then $I_{a,b,c}(f)$ maps $\mathcal{R}(\beta)$ into $ S_{p} $ class.
\ethm
\bpf Let $a \in {\Bbb C} \backslash \{ 0 \} ,\, b > 0,\, c > 0,\,   |a| \neq 1 ,\, b\neq 1$,\, $c \neq 1, \, b > |a|-1$ and $c > |a|-1$. \\

Consider the integral operator $\mathcal{I}_{a,b,c}(f)$  is given by (\ref{inteq7}). In the view of the equation (\ref{lem2eq1}), it is enough to show that
\begin{eqnarray*}
  T &=& \sum_{n=2}^{\infty}(2n-1)|B_n|\leq 1.
\end{eqnarray*}
(or)
\begin{eqnarray*}
  T &=& \sum_{n=2}^{\infty}(2n-1)\left|\frac{(a)_{n-1}(b)_{n-1}(c)_{n-1}}{(b+1)_{n-1}(c+1)_{n-1}(1)_{n-1}}\right|\, |a_n|\leq 1.
\end{eqnarray*}
Using the inequality $|(a)_n|\leq  (|a|)_n$ and the equation $(\ref{thm2eq2})$ in above, we have
\begin{eqnarray*}
  T  &\leq&  2(1-\beta)\left[2\sum_{n=0}^{\infty}  \frac{(|a|)_{n}(b)_{n}(c)_{n}}{(b+1)_{n}(c+1)_{n}(1)_{n}}-\sum_{n=0}^{\infty}  \frac{(|a|)_{n}(b)_{n}(c)_{n}}{(b+1)_{n}(c+1)_{n}(1)_{n+1}} -1  \right]:=T_1
\end{eqnarray*}
Using the equation (\ref{inteq6}) and the result (3) of Lemma \ref{lem1eq1}, we find that
\begin{eqnarray*}
  T_1 &\leq&  2(1-\beta)\left[ 2\frac{bc}{c-b}\Gamma(1-|a|)\left(\frac{\Gamma(b)}{\Gamma(1-|a|+b)}-\frac{\Gamma(c)}{\Gamma(1-|a|+c)}\right)  \right.\\ && \left.- \frac{bc}{(|a|-1)(b-1)(c-1)}\left(\left(\frac{\Gamma(2-a)}{c-b}\right)\left(\frac{(c-1)\Gamma(b)}{\Gamma(1-|a|+b)}-\frac{(b-1)\Gamma(c)}{\Gamma(1-|a|+c)}\right) -1 \right)-1  \right]
\end{eqnarray*}
Using the fact that $\Gamma(a+1)=a\Gamma(a)$, the above reduces to
\begin{eqnarray*}
 T &\leq&  2(1-\beta)\left[ \frac{bc}{c-b}\Gamma(1-|a|)\left(\frac{\Gamma(b)}{\Gamma(1-|a|+b)}\left( \frac{2b-1}{b-1}\right) \right.\right.\\ && \left.\left.- \frac{\Gamma(c)}{\Gamma(1-|a|+c)} \left(\frac{2c-1}{c-1}\right)\right) + \frac{bc}{(|a|-1)(b-1)(c-1)}-1\right]
\end{eqnarray*}
By the condition (\ref{thm5eq1}), the above expression is bounded by 1 and hence
\begin{eqnarray*}
 2(1-\beta)\left[ \frac{bc}{c-b}\Gamma(1-|a|)\left(\frac{\Gamma(b)}{\Gamma(1-|a|+b)}\left(\frac{2b-1}{b-1}\right)-\frac{\Gamma(c)}{\Gamma(1-|a|+c)} \left(\frac{2c-1}{c-1}\right)\right)\right.&&\\  \left.
+ \frac{bc}{(|a|-1)(b-1)(c-1)}-1\right] &\leq& 1
\end{eqnarray*}
Under the stated condition, the operator $\mathcal{I}_{a,b,c}(f)$ maps $\mathcal{R}(\beta)$ into $\es_p$ and  the proof is complete.
\epf
\bthm\label{thm6eq0}  Let $a \in {\Bbb C} \backslash \{ 0 \} ,\,  b > |a|-1$ and $c > |a|-1$. Suppose that $a,\, b$ and $c$ satisfy the condition
 \beq\label{thm6eq1}
 \frac{bc}{c-b}\, \Gamma(1-|a|)\left(\frac{(1-b)(1-2b)\Gamma(b)}{\Gamma(1-|a|+b)}-\frac{(1-c)(1-2c)\Gamma(c)}{\Gamma(1-|a|+c)}\right) &\leq& 2. \nonumber
\eeq
then $I_{a,b,c}(f)$ maps $\es$ into $ \es_{p} $ class.
\ethm
\bpf
The proof is similar to Theorem \ref{thm7eq1}. So we omit the details.
\epf
\section{Conclusion and Future Scope}
We derived the geometric properties for the Clausen's Hypergeometric Series $_3F_2(a,b,c;b+1,c+1;z)$ in which the numerator and denominator parameters differs by arbitrary negative integers.\\

It would be of great interest to determine conditions on the  parameters $a,b,c,\lambda $ and $\beta$ so that the  integral operator defined by (\ref{inteq7}) associated with some classes of univalent functions using Dixon's summation formula which we state as an open problem.\\

{\bf Problem:} To determine the conditions on the parameters $a,\,b,\,c,$ such that the hypergeometric function $_3F_2(a,b,c;1+a-b,\, 1+a-c)$ associated with the Dixon's summation formula or its equivalence has the admissibility to be the classes like $\mathcal{S}^{*}_{\lambda}, \, \mathcal{C}_{\lambda}$,\, $\mathcal{UCV}$  and $\es_p$.

\end{document}